\documentclass[12pt]{amsart}
\usepackage{amssymb,latexsym}

\newcommand{\ben}{\begin{enumerate}}
\newcommand{\een}{\end{enumerate}}
\newcommand{\ble}{\begin{lem}}
\newcommand{\ele}{\end{lem}}
\newcommand{\bth}{\begin{thm}}
\renewcommand{\eth}{\end{thm}}
\newcommand{\bpr}{\begin{prop}}
\newcommand{\epr}{\end{prop}}
\newcommand{\bco}{\begin{cor}}
\newcommand{\eco}{\end{cor}}
\newcommand{\bcon}{\begin{conj}}
\newcommand{\econ}{\end{conj}}
\newcommand{\bde}{\begin{defn}}
\newcommand{\ede}{\end{defn}}
\newcommand{\bex}{\begin{exa}}
\newcommand{\eex}{\end{exa}}
\newcommand{\barr}{\begin{array}}
\newcommand{\earr}{\end{array}}
\newcommand{\btab}{\begin{tabular}}
\newcommand{\etab}{\end{tabular}}
\newcommand{\beq}{\begin{equation}}
\newcommand{\eeq}{\end{equation}}
\newcommand{\bea}{\begin{eqnarray*}}
\newcommand{\eea}{\end{eqnarray*}}
\newcommand{\bce}{\begin{center}}
\newcommand{\ece}{\end{center}}
\newcommand{\bpi}{\begin{picture}}
\newcommand{\epi}{\end{picture}}
\newcommand{\bfi}{\begin{figure} \begin{center}}
\newcommand{\efi}{\end{center} \end{figure}}
\newcommand{\capt}{\caption}
\newcommand{\bsl}{\begin{slide}{}}
\newcommand{\esl}{\end{slide}}

\newcommand{\bib}{thebibliography}

\newcommand{\ol}{\overline}

\newcommand{\hs}[1]{\hspace{#1}}
\newcommand{\hso}[1]{\hspace{-1pt}}

\newcommand{\qmq}[1]{\quad\mbox{#1}\quad}

\newcommand{\sbs}{\subset}
\newcommand{\sbe}{\subseteq}

\newcommand{\Fh}{\hat{F}}

\newcommand{\lte}{\unlhd}

\newcommand{\jn}{\vee}

\newcommand{\mt}{\wedge}

\newcommand{\case}[4]{\left\{\barr{ll}#1&\mbox{#2}\\#3&\mbox{#4}\earr\right.}
\newcommand{\fl}[1]{\lfloor #1 \rfloor}

\def\<{\langle}
\def\>{\rangle}

\newcommand{\ree}[1]{(\ref{#1})}

\newcommand{\ra}{\rightarrow}

\newcommand{\de}{\delta}

\newcommand{\om}{\omega}

\newcommand{\Om}{\Omega}

\newcommand{\bbC}{{\mathbb C}}

\newcommand{\bbP}{{\mathbb P}}

\newcommand{\pb}{\ol{p}}
\newcommand{\Pb}{\ol{P}}

\newcommand{\aim}{Adv.\ in Math.\/}

\newcommand{\Gda}{\put(30,0){\circle*{3}}}

\newcommand{\Gag}{\put(0,60){\circle*{3}}}

\newcommand{\Ggg}{\put(60,60){\circle*{3}}}

\newcommand{\Gdaag}{\put(30,0){\line(-1,2){30}}}

\newcommand{\Gdagg}{\put(30,0){\line(1,2){30}}}

\newtheorem{thm}{Theorem}[section]
\newtheorem{prop}[thm]{Proposition}
\newtheorem{cor}[thm]{Corollary}
\newtheorem{lem}[thm]{Lemma}
\newtheorem{conj}[thm]{Conjecture}
\newtheorem{exa}[thm]{Example}

\newcommand{\Faa}{F(a,a)}
\newcommand{\Fab}{F(a,b)}
\newcommand{\Fac}{F(a,c)}
\newcommand{\Fbb}{F(b,b)}
\newcommand{\Fcc}{F(c,c)}
\newcommand{\blk}{\mathop{\rm bk}\nolimits}

\begin{document}
\title{GCD matrices, posets, and nonintersecting paths
}

\author{Ercan Altinisik}
\address{Department of Mathematics, Faculty of Arts and Sciences, Gazi
  University, 06500 Teknikokullar - Ankara, TURKEY} 
\email{ealtinisik@gazi.edu.tr}

\author{Bruce E. Sagan}
\address{Department of Mathematics, Michigan State University,
East Lansing, MI 48824-1027, USA}
\email{sagan@math.msu.edu}

\author{Naim Tuglu}
\address{Department of Mathematics, Faculty of Arts and Sciences, Gazi
  University, 06500 Teknikokullar - Ankara, TURKEY} 
\email{naimtuglu@gazi.edu.tr}

\date{\today}

	\begin{abstract}
We show that with any finite partially ordered set, $P$, one can
associate a matrix whose determinant factors nicely.  As corollaries,
we obtain a number of results in the literature about GCD matrices
and their relatives.  Our main theorem is proved combinatorially using
nonintersecting paths in a directed graph.
	\end{abstract}

\maketitle

\section{Introduction}
\label{i}

Let $\bbP$ denote the positive integers and suppose we are given a
subset $S=\{a_1,\ldots,a_n\}$ of $\bbP$.  The corresponding {\it GCD
matrix} is $(S)=(s_{ij})$ where $s_{ij}=(a_i,a_j)$, the greatest
common divisor of $a_i$
and $a_j$.  We say that $S$ is {\it factor closed\/} if given any
$a_i\in S$ and any divisor $d|a_i$ then $d\in S$.  
H. J. S. Smith~\cite{smi:vca}
proved the following beautiful result about the determinant of $(S)$.
\bth[Smith]
\label{smith}
If $S=\{a_1,\ldots,a_n\}$ is factor closed then
\beq
\label{(S)}
\det (S)=\phi(a_1)\cdots\phi(a_n)
\eeq
where $\phi$ is Euler's totient function.\qed
\eth

Since Smith's pioneering paper, a host of related results have
appeared in the literature.  For a survey with references, see the
paper of Haukkanen, Wang, and Sillanp\"a\"a~\cite{hwp:sd}.  We will
show that many of these are special cases of a general determinantal
identity associated with any finite partially ordered set $P$ (see
Theorem~\ref{main} below).  This result can be proved by factoring the matrix
constructed from $P$, but we choose to give a combinatorial proof
based on counting families of nonintersecting paths in acyclic digraphs.  This
technique is due to Lindstr\"om~\cite{lin:vri} and independently 
Karlin~\cite{kar:cpa}. It was later rediscovered and popularized by
Gessel and Viennot~\cite{gv:bdp,gv:dpp}.  We will use the approach in
Stembridge's paper~\cite{ste:npp}.

The rest of this article is structured as follows.  The next section
will be devoted to  proving our main theorem and giving some
preliminary applications.  In Section~\ref{ms}, we specialize to the
case where $P$ is a meet semilattice.  This permits us to invert the
sums appearing in the general case.  Section~\ref{lcmcs} investigates
what happens to the determinant if the set $S$ is lower closed (the
partially ordered set generalization of being factor closed) or meet closed.  
We end with a section containing comments and an open question.

\section{The main theorem}
\label{mt}

We will first review the digraph machinery which we will need to prove
our main theorem.  Explanations of any undefined terms from graph
theory can be
found in the texts of Harary~\cite{har:gt} or Chartrand and
Lesniak~\cite{cl:gd}.  For definitions of concepts about
posets (partially ordered sets), the reader can consult Stanley's
book~\cite{sta:ec1}. 

Let $D$ be a finite acyclic digraph with vertices $V$ and arcs $A$.
Let $R$ be a commutative ring with identity and suppose we are given a  function (weighting)
$\om:A\ra R$.  
Then any directed path $p:\ v_0v_1\ldots v_k$ has an associated weight
$$
\om(p)=\prod_{i=1}^k \om(v_{i-1}v_i).
$$  
We also let 
$$
\om(v_0,v_k)=\sum_p \om(p)
$$
where the sum is over all directed paths $p$ from $v_0$ to $v_k$.  Note
that this sum is finite because $D$ is finite and acyclic.  

Now suppose we are given two disjoint sets of vertices
$V'=\{v_1',\ldots,v_n'\}$ and $V''=\{v_1'',\ldots,v_n''\}$ in $V$.  Consider
an $n$-tuple of directed paths $\pi=(p_1,\ldots,p_n)$ where $p_i$ goes
from $v_i'$ to $v_i''$ for all $i$, $1\le i\le n$.  Then we assign
weights to $\pi$ and to the pair $(V',V'')$ in a way analogous to the
one used in the preceding paragraph
$$
\om(\pi)=\prod_{i=1}^n \om(p_i)\qmq{and} \om(V',V'')=\sum_{\pi}\om(\pi)
$$ 
where the sum is over all $\pi$ where no two of the paths in the
$n$-tuple intersect.  

Finally, given a permutation $g$ in the symmetric group $S_n$, we
define $\pi_g=(p_1,\ldots,p_n)$ to be an $n$-tuple of directed paths
such that $p_i$ goes from $v_i'$ to $v_{g(i)}''$ for all $i$.  So the
$n$-tuples considered in the previous paragraph would have $g=e$
where $e$  is the identity permutation.  If any pair of 
paths in $\pi_g$ intersect then we call the $n$-tuple {\it
intersecting\/}  and {\it nonintersecting\/} otherwise.
We can now state Stembridge's theorem~\cite{ste:npp} about enumerating
nonintersecting paths.

\bth [Stembridge] 
\label{paths}
Let $D$ be a finite acyclic digraph with disjoint vertex sets 
$V'=\{v_1',\ldots,v_n'\}$ and $V''=\{v_1'',\ldots,v_n''\}$ such that
any $\pi_g$ is intersecting if $g\neq e$.   Let
$(D)=(d_{ij})$ be the matrix with
$$
d_{ij}=\om(v_i',v_j'').
$$
Then
$$
\det (D)=\om(V',V'').\qed
$$
\eth

Now let $P$ be a finite poset and consider the {\it incidence
algebra\/} $I(P,R)$ of $P$ over $R$ which consists of all functions
$F:P\times P\ra R$ such that $F(a,b)=0$ unless $a\le b$.  The identity
element of $I(P,R)$ is the Kronecker delta function 
$$
\de(a,b)=\case{1}{if $a=b$,}{0}{otherwise.}
$$ 
We will also need the {\it zeta function\/} of $I(P,R)$ defined by
$$
\zeta(a,b)=\case{1}{if $a\le b$,}{0}{otherwise.}
$$
The zeta function is invertible and its inverse is called the M\"obius
function $\mu$ of $I(P,R)$.  In other words, $\mu$ is the unique
function in $I(P,R)$ satisfying
$$
\sum_{a\le c\le b} \mu(a,c)=\sum_{a\le c\le b} \mu(c,b)=\de(a,b).
$$

Now fix a linear ordering of $P$ which can be used to index the rows and
columns of a matrix $M$.  Since we will be taking determinants, it will
not matter which linear order is used and we will merely say that the
matrix is indexed by $P$.  
We can now state and prove our main theorem.
\bth
\label{main}
Let $P$ be a finite poset and let $F,G\in I(P,R)$.
Let $(P)_{FG}$ be the  matrix indexed by $P$ with entries
\beq
\label{pab}
p_{ab}=\sum_{c\in P} F(c,a)G(c,b).
\eeq
Then
$$
\det (P)_{FG}=\prod_{a\in P} F(a,a)G(a,a).
$$
\eth
\proof
Construct a digraph $D$ as follows.  For the vertices of $D$, take
three copies of the elements of $P$ which we denote by $P'$, $P''$,
and $P'''$.  
Now put an arc from $a'\in P'$  to $c'''\in P'''$ if and only if $a\ge c$ in
$P$.  Dually, put an arc from $c'''\in P'''$ to $b''\in P''$ if and only if
$c\le b$.  Finally, give weights to the arcs by
$$
\om(a',c''')=F(c,a) \qmq{and} \om(c''',b'')=G(c,b).
$$

Consider paths $p$ from $a'\in P'$ to $b''\in P''$.  Clearly all such
paths have the form $p: a',c''',b''$ where $c\le a$ and $c\le b$ in
$P$.  So if we take $V'=P'$ and $V''=P''$ then we have $(D)=(P)_{FG}$ since
$$
d_{a'b''}=\om(a',b'')
=\sum_{c\le a,\ c\le b } F(c,a)G(c,b) =p_{ab}.
$$
The last step follows from the fact that a term in $p_{ab}$
corresponding to some $c$ not satisfying the given inequalities is zero.

We need to show that $D$ satisfies the hypotheses of
Theorem~\ref{paths}.  By construction, $D$ is acyclic.  Now
take any $n$-tuple of paths $\pi_g$ which is 
nonintersecting.  We will show that this forces $g=e$.  First we claim
that the path starting at any $a'\in P'$ must then follow the arc to
$a'''\in P'''$.  Suppose this is not true for some $a'$ and let
$c'''$ be the next vertex on the path.  Then
we must have $a>c$ in $P$.  Now consider the path of $\pi_g$ starting
at $c'$.  Since it cannot intersect the previous path, it must go to
some $d'''$ with $c>d$.  Continuing in this fashion, we can construct
an infinite decreasing chain in $P$, contradicting the assumption that
$P$ is finite.  So our claim is true.  By a dual argument, one can show
that the path starting at $a'$ must continue from $a'''$ to $a''$ and
so $g=e$.  Furthermore, we have shown that this family $\pi_e$ is the only
nonintersecting path family.   So, by the way we have defined the weights and
Theorem~\ref{paths}, 
$$
\det (D) = \om(\pi_e)=\prod_{a\in P} F(a,a)G(a,a).\qed
$$

We should note that this theorem can also be proved
algebraically, either as a corollary to a result of
D. A. Smith~\cite[Corollary 2]{smi:bfa} about the incidence algebra of
$P$, or by directly factoring the matrix for $P$ as $(P)_{FG}=M_F^t M_G$ where
$M_F$ and $M_G$ are the matrices corresponding to the incidence
algebra elements $F$ and $G$, respectively, and $t$ denotes transpose.

\thicklines
\setlength{\unitlength}{1pt}
\bfi
\bpi(50,80)(0,-10)
\Gda \Gag \Ggg
\Gdaag  \Gdagg
\put(30,-10){\makebox(0,0){$a$}}
\put(0,70){\makebox(0,0){$b$}}
\put(60,70){\makebox(0,0){$c$}}
\epi
\capt{A poset $P$}\label{poset}
\efi

As an example of Theorem~\ref{main}, consider the poset $P$ whose
Hasse diagram is shown in Figure~\ref{poset}.  Then using the linear
order $a,b,c$ one obtains
$$
\barr{l}
\hs{-2pt}\det
  \left(
  \barr{ccc}
  \Faa G(a,a) &\Faa G(a,b)             &\Faa G(a,c)\\[5pt]
  \Fab G(a,a) &\Fab G(a,b)+\Fbb G(b,b) &\Fab G(a,c)\\[5pt]
  \Fac G(a,a) &\Fac G(a,b)             &\Fac G(a,c)+\Fcc G(c,c)
  \earr
  \right)\\[30pt]
\qquad=F(a,a)G(a,a)F(b,b)G(b,b)F(c,c)G(c,c)
\earr
$$

As a first application, we note a corollary of Theorem~\ref{main}
which simultaneously generalizes a theorem of Apostol~\cite{apo:apg}
and one of Daniloff~\cite{dan:ctf}.  Let $f,g$ be arbitrary functions
from the poset $P$ to the ring $R$.  Then on substituting $F(a,b)f(a)$
and $G(a,b)g(a)$ for $F(a,b)$ and $G(a,b)$, respectively, we
immediately have the following result.
\bco
\label{FfGg}
Let $P$ be a finite poset.  Let $(\Pb)$ be the matrix 
indexed by $P$ with entries
$$
\pb_{ab}=\sum_{c\in P} F(c,a)f(c)G(c,b)g(c).
$$
Then
$$
\det (\Pb)=\prod_{a\in P} F(a,a)f(a)G(a,a)g(a).\qed
$$
\eco

Now consider the poset defined by using the divisor ordering on
$P_n\stackrel{\rm def}{=}\{1,2,\ldots,n\}$.
Apostol's theorem is obtained by specializing the previous
result to the case where $P=P_n$,
$F(a,b)=\zeta(a,b)$, $G(a,b)=G(b/a)$ for any function $G:P_n\ra R$, and
$g(a)=1$ for all $a\in P$.  In this case the determinant becomes
$$
\det (\Pb_n)=f(1)f(2)\cdots f(n) G(1)^n.
$$
With this evaluation in hand, Apostol showed that letting $f(a)=a$ for
$a\in P_n$  and $G(a,b)=\mu(b/a)$ where $\mu$ is the usual
number-theoretic M\"obius function gives
$$
\det (c(a,b)) = n!
$$
where $c(a,b)$ is Ramanujan's sum.

To obtain Daniloff's theorem, let $P=P_n$, $g(a)=1$ for all $a\in P_n$,
and $F(a,b)=G(a,b)=\Om_k(b/a)$  where
$$
\Om_k(a)=\case{a^{1/k}}{if $a^{1/k}\in \bbP$,}{0}{else.}
$$
Now Corollary~\ref{FfGg} becomes Daniloff's result,
$$
\det (\Pb_n)=f(1)f(2)\cdots f(n). 
$$
Note that one obtains  the same evaluation not just for $\Om_k$, but also
for any functions $F,G:P\ra R$ such that $F(a)=G(a)=1$ for all $a\in P$.

\section{Meet semilattices}
\label{ms}

Now suppose that our poset is a meet semilattice $L$ so that every pair of
elements $a,b\in L$ have a greatest lower bound or {\it meet\/} 
$a\mt b$.  Note that in this case the sum~\ree{pab} can be restricted
to $c\le a\mt b$.  This special case of Theorem~\ref{main} was
discovered by Haukkanen, Wang, and Sillanp\"a\"a~\cite{hwp:sd}.  They
also showed that by further specialization one could obtain a theorem
of Jager~\cite{jag:uai} concerning a unitary analogue of~\ree{(S)},
Smith's evaluation~\cite{smi:vca} of the LCM determinant, as well as
a number of other results.  

It would be nice to be able to compute the value of 
determinants where each entry is a single term, as in Smith's original
case, rather than a sum.
Since~\ree{pab} now consists of a sum over all $c$ below a certain
element in $L$, one can use M\"obius inversion to accomplish this.
Thus we can prove
the following theorem of Lindstr\"om~\cite{lin:ds} as a corollary to
Theorem~\ref{main}. 
\bth[Lindstr\"om]
\label{main2}
Let $L$ be a finite meet semilattice and let $f\in I(L,R)$. 
Let $(L)_f$ be the matrix indexed by $L$ with entries
$$
l_{ab}=f(a\mt b,a)
$$
Then
\beq
\label{Lf}
\det (L)_f=\prod_{a\in L} \left(\sum_{c\in L}\mu(c,a)f(c,a)\right).
\eeq
\eth
\proof
Define a function $F\in I(L,R)$ by
$$
F(a,b)=\sum_{c\le a} \mu(c,a)f(c,b)
$$
if $a\le b$ and $F(a,b)=0$ otherwise.  By M\"obius inversion, this is
equivalent to
$$
f(a,b)=\sum_{c\le a} F(c,b).
$$
It follows that
$$
l_{ab}=f(a\mt b,a)=\sum_{c\le a\mt b} F(c,a)
=\sum_{c\in L} F(c,a)\zeta(c,b).
$$
So applying Theorem~\ref{main} we obtain
$$
\det (L)_f=\prod_{a\in L} F(a,a)\zeta(a,a)=
\prod_{a\in L} \left(\sum_{c\in L}\mu(c,a)f(c,a)\right). \qed
$$

This result can also be found as an exercise in Stanley's 
book~\cite[Chapter 3, Exercise 37]{sta:ec1}.  The special 
case where $f(a,b)$ depends only on $a$ was proved independently by
Wilf~\cite{wil:hdm}.  Smith himself~\cite{smi:vca} noted that this
theorem holds under the further assumption that $L$ is a
factor closed subset of $\bbP$.

To see how Smith's determinant~\ree{(S)} follows from
Theorem~\ref{main2}, just let $f(a,b)=a$ for all $a\in S$.  Then 
$s_{ij}=f(a_i\mt a_j, a_i)$ and, since $S$ is factor closed,
$$
\sum_{c\in S}\mu(c,a)f(c,a)=
\sum_{c|a}\mu(a/c)\ c =\phi(a).
$$

\section{Lower-closed and meet-closed sets}
\label{lcmcs}

In work on analogues of Theorem~\ref{smith}, there are two conditions
that are commonly imposed on the set $S$.  If $S\sbe P$ for some poset
$P$, we say that $S$ is {\it lower closed\/} or a {\it lower order
ideal\/} if $a\in S$ and $b\le a$ in $P$ implies $b\in S$.  This
corresponds exactly to $S$ being factor closed if it is a subset of
$\bbP$ ordered by division.  If, in particular, our poset is a
meet semilattice $L$ and $S$ is lower closed, then $S$ is also a meet
semilattice with the same M\"obius function as $L$.  So our previous
results cover this case without change.

If $S\sbe L$ for a meet semilattice one can also talk about $S$ being
{\it meet closed\/} which means that if $a,b\in S$ then $a\mt b\in S$
where the meet is taken in $L$.  Again, $S$ is also a meet semilattice
and so Theorem~\ref{main2} still applies and in fact generalizes a
result of Bhat~\cite{bha:gcd} who considered the case when $f$ is a
function of only one argument.

However, if $S\sbe L$ is meet closed we will show that there is also an
analogue of~\ree{Lf} which involves the M\"obius function of $L$,
rather than that of $S$.  This result will generalize a theorem of
Haukkanen~\cite{hau:mmp} who, as in Bhat's paper, only looked at
functions of a 
single variable.  Specializing yet further to $L=P_n$ and the function
$f(a)=a$ for all $a\in P_n$, one obtains a
theorem of Beslin and Ligh~\cite{bl:ags}.  

To state our result, it will be convenient to have some notation.
Let $\ell$ be a linear extension of the partial order on $S$
so that if  $\ell=a_1,a_2,\ldots,a_n$ then $a_i<a_j$ implies $i<j$.
If $d\in L$ then we will write
$d\lte a_i$ if $d\le a_i$ and $d\not\le a_j$ for any $j<i$.
\bth
Let $L$ be a finite meet semilattice and let $f\in I(L,R)$.  Suppose
$S\sbe L$ is meet closed and fix a linear extension
$\ell=a_1,a_2,\ldots,a_n$ of $S$.  Then 
$$
\det(S)_f=\prod_{i=1}^n 
\left(\sum_{d\lte a_i}\ \sum_{c\in L}\mu(c,d)f(c,a_i)\right)
$$
\eth
\proof
We will use unsubscripted variables for elements of $L$ which are not
necessarily in $S$.  Define $F\in I(L,R)$ as in the proof of
Theorem~\ref{main2} so that we have
\beq
\label{f}
f(a_i,a_j)=\sum_{d\le a_i} F(d,a_j).
\eeq
Also define $\Fh\in I(S,R)$ by
$$
\Fh(a_i,a_j)=\sum_{d\lte a_i}\ \sum_{c\in L}\mu(c,d) f(c,a_j)
=\sum_{d\lte a_i} F(d,a_j)
$$
if $a_i\le a_j$ and $\Fh(a_i,a_j)=0$  otherwise.  If we can show that
$$
f(a_i,a_j)=\sum_{a_k\le a_i} \Fh(a_k,a_j)
$$
then the rest of the proof will follow as in the demonstration of
Theorem~\ref{main2}.   So, by the definition of $\Fh$, it suffices to
show
\beq
\label{f2}
f(a_i,a_j)=\sum_{a_k\le a_i}\ \sum_{d\lte a_k} F(d,a_j).
\eeq
We will do this by showing that there is a one-to-one correspondence
between the terms in~\ree{f2} and those in~\ree{f}.

First note that each $d\in L$ occurs at most once in~\ree{f} and at
most once in~\ree{f2} (since $d\lte a_k$ for at most one $a_k\in S$).
If $d$ occurs in~\ree{f2} then $d\le a_k\le a_i$ and so $d$ occurs
in~\ree{f}.  Conversely, if $d$ occurs in~\ree{f} then $d\le a_i$ and
so we must have $d\lte a_k$ for some $a_k$ with $k\le i$.  But now
$d\le a_i\mt a_k = a_l$ for some $l$ since $S$ is meet closed.
Furthermore, $a_l\le a_k$ implies $l\le k$, and so $l=k$ since
$d\lte a_k$.  It follows that $d\lte a_k$ where $a_k=a_l\le a_i$.
Thus $d$ occurs in~\ree{f2} and we have finished the proof.\qed

We should note that Breslin and Ligh~\cite{bl:gcd} have derived a
formula for $(S)_f$ for any subset $S\sbs\bbP$ in terms of an arbitrary
lower-closed set containing $S$.  Since this identity expresses
$(S)_f$ as a product of two (not necessarily square) matrices, one can
appy the Cauchy-Binet formula to it, as done in Li's
paper~\cite{li:dgc}, to obtain $\det (S)_f$ as a sum of determinants.
This approach was also generalized to meet semi-lattices in an article
of Haukkanen~\cite{hau:mmp}.

\section{Comments and an open question}

Our results can also yield other
information about the matrices $(P)_{FG}$ and $(L)_f$.  The following
theorem is an example of this.  In it, we specialize the ring to be
the complex numbers $\bbC$ for simplicity, although more general rings
can be used.

\bth
Let $P$ be a finite poset and let $F,G\in I(P,\bbC)$.
\ben
\item  The matrix $(P)_{FG}$ is invertible if and only if
  $F(a,a),G(a,a)\neq0$ for all $a\in P$.
\item  The matrix $(P)_{FG}$ is positive definite if and only if
  $F(a,a)G(a,a)>0$ for all $a\in P$.
\een
\eth
\proof
Part (1) follows immediately from Theorem~\ref{main} and the fact that
a matrix over $\bbC$ is invertible if and only if its determinant is
nonzero.

For part (2), let the total order used to index the rows and columns
of $P$ be a linear extension of $P$.  Then each principal
submatrix of $P$ is indexed by a lower-closed subset $S$ of $P$.  It
follows that the submatrix indexed by $S$ is exactly $(S)_{FG}$.  But
now we are done by Theorem~\ref{main} again, since a matrix is
positive definite if and only if the principal subdeterminants are all
positive.  \qed 

If $L$ is a finite lattice, then one can compute the
determinant of a matrix involving joins (least upper bounds) $a\jn b$
of elements $a,b\in L$ by working in the dual of $L$ so that joins
become meets.  
For general results about join and meet matrices, see the paper of
Korkee and Haukkanen~\cite{kh:mjm}.

However, there are factorizations of such determinants which
we have been unable to obtain by our methods.  As an example, we
consider the matrix $T_n(q)$ of chromatic joins introduced by Birkhoff
and Lewis~\cite{bl:cp} to study the chromatic polynomial of planar maps.
To define this matrix, let $q$ be a formal parameter. 
Let $\Pi_n$ denote the lattice of partitions of $\{1,2,\ldots,n\}$
ordered by refinement and let $NC_n$ denote the lattice of noncrossing
partitions of the same set.  (An excellent survey about noncrossing
partitions can be found in the article of Simion~\cite{sim:np}.)
Then Tutte~\cite{tut:ble} showed that the matrix of chromatic joins
could be defined as
$$
T_n(q)=\left(q^{\blk(a\jn_{\Pi_n} b)}\right)_{a,b\in NC_n}
$$
where $\blk(a)$ is the number of blocks (subsets) in the partition $a$.

Tutte~\cite{tut:mcj} also derived a product formula for the determinant
of $T_n(q)$ in terms of Beraha polynomials.  (See also
Dahab~\cite{dah:ble}.)  Letting $\fl{\cdot}$ denote the floor or round
down function, the {\it $n$th Beraha polynomial\/} is defined to be
$p_0(q)=0$ for $n=0$, and for $n\ge1$
$$
p_n(q) =\sum_{i=0}^{\fl{n/2}} (-1)^i{n-i-1\choose i}q^{\fl{n/2}-i}.
$$
Using the version of Tutte's formula in the paper of Copeland,
Schmidt, and Simion~\cite{css:ntd} gives
$$
\det T_n(q) = q^{{2n-1\choose n}}   \prod_{m=1}^{n-1}
\left[\frac{p_{m+2}(q)}{qp_m(q)}\right]^{\frac{m+1}{n}{2n\choose n-m-1}}.
$$
The same paper also contains a related determinant-product identity
which the authors 
note that they were unable to prove using Lindstr\"om's Theorem
(Theorem~\ref{main2} above), although other proofs exist.  It would be
interesting to find a way to apply the machinery of this paper to
these identities.

\medskip

{\it Acknowledgement.}  We would like to thank Vic Reiner for helpful
discussions and, in particular, referring us to the paper of Copeland,
Schmidt, and Simion~\cite{css:ntd}.

\begin{\bib}{99}

\bibitem{apo:apg} T. M. Apostol, Arithmetical properties of
generalized Ramanujan sums, {\it Pacific J. Math.\ } {\bf 41} (1972),
281--293. 

\bibitem{bl:ags} S. Beslin and S. Ligh, Another generalization of
Smith's determinant, {\it Bull.\ Austral.\ Math.\ Soc.\ } {\bf 40}
(1989), 413--415.

\bibitem{bl:gcd} S. Beslin and S. Ligh, Greatest common divisor
matrices, {\it Linear Algebra Appl.\/} {\bf 118}
(1989), 69--76.

\bibitem{bha:gcd} B. V. Bhat, On greatest common divisor matrices and
their applications, {\it Linear Algebra Appl.\/} {\bf 158} (1991),
77--97.  

\bibitem{bl:cp} G. D. Birkhoff and D. C. Lewis, Chromatic polynomials,
{\it Trans.\ Amer.\ Math.\ Soc.\ } {\bf 60} (1946), 355--451.

\bibitem{cl:gd} G. Chartrand and L. Lesniak, ``Graphs and Digraphs,'' 
second edition, Wadsworth \& Brooks/Cole, Pacific Grove, CA, 1986.

\bibitem{css:ntd} A. Copeland, F. Schmidt, and R. Simion, Note on two
determinants with interesting factorizations,
{\it Discrete Math.\/} {\bf 256} (2002), 449--458.

\bibitem{dah:ble} R. Dahab, ``The Birkhoff-Lewis Equations,'' Ph.D. thesis,
University of Waterloo, Waterloo, ONT, 1993. 

\bibitem{dan:ctf} G. Daniloff, Contribution \`a la th\'eorie des
fonctions arithm\'etiques (in Bulgarian; French summary), {\it Sb.\
Bulgar.\ Akad.\ Nauk\/} {\bf 35} (1942), 479--590.

\bibitem{gv:bdp} I.  Gessel  and  G.  Viennot,  Binomial 
determinants, paths, and hook length formulae, {\it \aim} {\bf 58}
(1985), 300--321.

\bibitem{gv:dpp} I. Gessel and G. Viennot, Determinants, paths, and
plane partitions, in preparation.

\bibitem{har:gt} F. Harary, ``Graph Theory,'' Addison-Wesley, Reading,
MA, 1971.

\bibitem{hau:mmp} P. Haukkanen, On meet matrices on posets,
{\it Linear Algebra Appl.\/} {\bf 249} (1996), 111--123.  

\bibitem{hwp:sd} P. Haukkanen, J. Wang, and J. Sillanp\"a\"a, On
Smith's determinant, {\it Linear Algebra Appl.\/} {\bf 258} (1997),
251--269.  

\bibitem{jag:uai} H. Jager, The unitary analogues of some identities
for certain arithmetical functions, {\it Neder.\ Akad.\ Wetensch.\
Proc.\ Ser.\ A\/} {\bf 64} (1961), 508--515.

\bibitem{kar:cpa} S. Karlin, Coincident probabilities and applications
to combinatorics, {\it J. Appl.\  Probab.\ } {\bf 25A} (1988), 185--200.

\bibitem{kh:mjm} I. Korkee and P. Haukkanen, On meet matrices and join
matrices associated with incidence functions, {\it Linear Algebra
Appl.\/} {\bf 372} (2003), 127--153.

\bibitem{li:dgc} Z. Li, The determinants of GCD matrices, {\it Linear
Algebra Appl.\/} {\bf 134} (1990), 137--143.

\bibitem{lin:ds} B. Lindstr\"om, Determinants on semilattices, 
{\it Proc.\ Amer.\ Math.\ Soc.\ } {\bf 20} (1969), 207--208.

\bibitem{lin:vri} B. Lindstr\"om, On the vector representation of induced
matroids, {\it Bull.\  London Math.\  Soc.\ } {\bf 5} (1973), 85--90.

\bibitem{rot:tmf} G.-C. Rota, On the foundations of combinatorial
theory I. Theory of M\"obius functions, {\it Z. 
Wahrscheinlichkeitstheorie} {\bf 2} (1964), 340--368.

\bibitem{sim:np} R. Simion, Noncrossing partitions,
{\it Discrete Math.\ } {\bf 217} (2000), 367--409.

\bibitem{smi:bfa} D. A. Smith, Bivariate function algebras on posets,
{\it J. Reine Angew.\ Math.\/} {\bf 251} (1971) , 100--109.

\bibitem{smi:vca} H. J. S. Smith, On the value of a certain
arithmetical determinant, {\it Proc.\ London Math.\ Soc. Ser. 1\/}  
{\bf 7} (1876), 208--212.

\bibitem{sta:ec1} R. P. Stanley, ``Enumerative Combinatorics,
Volume 1,''  Cambridge University Press, Cambridge, 1997.

\bibitem{ste:npp} J. R. Stembridge, Nonintersecting paths, pfaffians
and plane partitions, {\it Adv.\ in Math.\ } {\bf 83} (1990), 96--131.

\bibitem{tut:mcj} W. T. Tutte, The matrix of chromatic joins,
{\it J. Combin.\ Theory Ser.\ B\/} {\bf 57} (1993), 269--288.

\bibitem{tut:ble} W. T. Tutte, On the Birkhoff-Lewis equation,
{\it Discrete Math.\ } {\bf 92} (1991), 417--425.

\bibitem{wil:hdm} H. S. Wilf, Hadamard determinants, m\"obius
functions, and the chromatic number of a graph, {\it Bull.\ Amer.\
Math.\ Soc.\/} {\bf 74} (1968), 960--964.

\end{\bib}

\end{document}